\documentclass{article}
\usepackage{url}
\usepackage{amsthm}
\usepackage{amsfonts}
\usepackage{amssymb}
\usepackage{amsgen}
\usepackage{amsmath}
\usepackage{amsopn}
\usepackage{verbatim}
\usepackage{xypic}
\usepackage{xspace}
\usepackage{multicol}
\usepackage{makeidx}
\usepackage{eepic}
\usepackage{upref}

\theoremstyle{plain}
\newtheorem{thm}{Theorem~}[section]
\newtheorem{lem}[thm]{Lemma}
\newtheorem{prop}[thm]{Proposition}
\newtheorem{cor}[thm]{Corollary}
\newtheorem{conj}[thm]{Conjecture}
\theoremstyle{definition}

\newtheorem{defn}[thm]{Definition}
\newtheorem{rem}[thm]{Remark}

\newcommand{\tn}{{\tilde{n}}}

\newcommand{\ors}{{\bold{s}}}

\newcommand{\N}{{\mathbb N}}
\newcommand{\ZN}{{{\mathbb Z}_{\ge 0}}}
\newcommand{\Q}{{\mathbb Q}}

\DeclareMathOperator{\pari}{{\mathfrak p}}

\newcommand{\ga}{{\alpha}}
\newcommand{\gl}{{\lambda}}
\newcommand{\gb}{{\beta}}
\newcommand{\gd}{{\delta}}

\newcommand{\gt}{{\tau}}

\begin{document}
\title{\bf  Finiteness of $p$-Divisible Sets of \\
Multiple Harmonic Sums}
\author{Jianqiang Zhao}
\date{}
\maketitle

\medskip\noindent\textbf{Abstract.}
Let $l$ be a positive integer and $\ors=(s_1,\dots,s_l)$ be a
sequence of positive integers. In this paper we shall study the
arithmetic properties of multiple harmonic sum $H(\ors;n)$ which
is the $n$-th partial sum of multiple zeta value series
$\zeta(\ors)$. We conjecture that for every $\ors$ and every prime
$p$ there are only finitely many $p$-integral partial sums
$H(\ors;n)$. This generalizes a conjecture of Eswarathasan and
Levine and Boyd for harmonic series. We provide a lot of evidence
for this general conjecture as well as some heuristic argument to
support it. This paper is a sequel to
\emph{Wolstenholme Type Theorem for multiple harmonic
sums}, Intl.\ J.\ of Num.\ Thy.\ \textbf{4}(1) (2008) 73-106.

\medskip\noindent\textbf{R\'esum\'e.}
Soit $l$ un entier et $\ors=(s_1, \dots, s_l)$ une
s\'equence d'entiers positifs. Dans ce document, nous \'etudierons les
propri\'et\'es arithm\'etique de sommes harmoniques multiples
$H(\ors; n)$, qui est le $n$-ème somme partielle de la valeur de
la s\'erie multiple zeta $\zeta(\ors)$. On conjecture que pour tout
$\ors$ et de tous les premiers
$p$, il n'y a que de nombreux finitely $p$-partie int\'egrante sommes
$H(\ors,n)$. Ceci g\'en\'eralise une conjecture de Eswarathasan et
Levine et Boyd pour la s\'erie harmonique. Nous fournissons beaucoup
d'\'el\'ements de preuve pour cette conjecture g\'en\'erale
ainsi que certaines heuristiques argument soutenir. Ce document fait suite \`a
\emph{Wolstenholme Type Theorem for multiple harmonic
sums}, Intl.\ J.\ of Num.\ Thy.\ \textbf{4}(1) (2008) 73-106.



\section{Introduction}
In recent years there is a revival of interest in the multiple zeta values
defined by
\begin{equation*}
\zeta(\ors):=\zeta(s_1,\dots, s_l)=\sum_{0<k_ 1<\dots<k_l}
k_1^{-s_1}\cdots k_l^{-s_l}
\end{equation*}
for $\ors=(s_1,\dots, s_l)\in \N^l$, where $\N$ is the set of
positive integers (see, for e.g., \cite{H,IKZ,Zag}).
In physics, however, not only these series but
also their partial sums have significant meanings in applications
(see \cite{Blum1,Blum2}). These partial sums are called the multiple
harmonic sums (MHS for short) which generalize the notion
of harmonic sums. In general it is defined as
\begin{equation}\label{equ:def}
H(\ors;n):=\sum_{1\le k_ 1<\dots<k_l\le n} k_1^{-s_1}\cdots
k_l^{-s_l},\quad n\in \ZN.
\end{equation}
where $|\ors|:=s_1+\dots+s_l$ is called the weight and $l=l(\ors)$ the length.
By convention we set $H(\ors;r)=0$ for $r=0,\dots,l-1$, and
$H(\emptyset;0)=1$. To save space, we denote by $\{e_1,\dots, e_t\}^d$
the string formed by repeating $(e_1,\dots,e_t)$ $d$ times and
$1^d=\{1\}^d$. We also write $H(\ors)$ for the set of all the
partial sums $H(\ors;n)$ when no confusion is likely to arise .

This paper can be considered as a sequel to \cite{1stpart} whose main goal
is to provide generalizations of Wolstenholme's Theorem for the MHS.
Hoffman \cite{Hmodp} obtained similar types of results independently.
One of the results in \cite{1stpart} is the following generalization to
homogeneous MHS:
\begin{thm} \label{th:p-1} {\em \cite[Theorem~2.13]{1stpart})}
Let $s$ and $l$ be two positive integers.  Let $p$ be an odd prime
such that $p\ge l+2$ and $p-1$ divides none of $ks$ and $ks+1$ for
$k=1,\dots, l$. Then
$$H(\{s\}^l; p-1) \equiv \begin{cases}
0 \quad &\pmod{\,p\,} \text{ if $ls$ is even,}\\
0 \quad &\pmod{p^2} \text{ if $ls$ is odd}.
\end{cases}$$
In particular, the above is always true if $p\ge ls+3$.
\end{thm}

One can also investigate the sums $H(\ors; n)$ with fixed $\ors$
but varying $n$. Such a study for harmonic series was initiated
systematically by Eswarathasan and Levine \cite{EL} and Boyd
\cite{Boyd}, independently. It turns out that to obtain precise
information one has to study Wolstenholme type congruences in some
detail and so these two directions of research are interwoven into
each other rather tightly. To state our main results and
conjectures we define
$$H(\ors;n)= \frac{a(\ors;n)}{b(\ors;n)},
\quad a(\ors;n), \ b(\ors;n)\in\N, \  \ {\rm
gcd}(a(\ors;n),b(\ors;n))=1.$$ For completeness, we set
$a(\ors;0)=0$ and $b(\ors;0)=1.$ Fixing a prime $p$ we are
interested in the $p$-divisibility of the integers $a(\ors;n)$ and
$b(\ors;n)$ for varying $n$. Thus we put $H(\ors;n)$ inside
$\Q_p$, the fractional field of the $p$-adic integers and let
$v_p$ be the discrete valuation on $\Q_p$ such that $v_p(p)=1$. In
this general situation we're forced to change the notation used by
the previous authors. For any $m\in \N$ and $\ors\in \N^l$, put
\begin{align*}
I(\ors|m):=&\{n\in \ZN: b(\ors;n)\not\equiv 0\pmod m\},\\
J(\ors|m):=&\{n\in \ZN: a(\ors;n)\equiv 0\pmod m\}.
\end{align*}
Note that $J(\ors|m)\ne \emptyset$ since $0\in J(\ors|m)$ always.
For any prime $p$ we call $J(\ors|p)$ the {\em $p$-divisible set}
of the MHS $H(\ors)$ defined by Eq.~\eqref{equ:def}.

In \cite{Boyd} Boyd presented a heuristic argument by modeling on
simple branching processes to convince us that the $p$-divisible
set of the harmonic series is finite for every prime $p$ (this is
also independently conjectured by Eswarathasan and Levine
\cite[Conjecture A]{EL}). Boyd also proves this conjecture for all
primes less than 550 except  for 83, 127 and 397. We now provide a generalization:
\begin{conj} \label{conj:main}
For any $\ors$ and any prime $p$ the $p$-divisible set $J(\ors|p)$
is finite.
\end{conj}
Although we are not able to prove this conjecture in general, we
obtain a lot of partial results. The primary tool to prove these
when $l(\ors)\ge 2$ is our Criterion Theorem~\ref{thm:Jfinite}.
Fixing an arbitrary prime $p$ we define
$$G_0=\{0\}  \text{ and } G_t=\{n: p^{t-1}\le n< p^t\} \text{ for } t\in\N.$$

\noindent\textit{Criterion Theorem.}  {\em Let $l\ge 2$ be a
positive integer and $p$ be a prime such that $l\in G_{t_0}$. Let
$\ors=(s_1,\cdots,s_l)\in \N^l$ and put $m=\min\{s_i: 1\le i\le
l\}$. For $t\in \N$ set
$f(\ors,p;t)=\min\{-v_p\big(H(\ors;n)\big): n\in G_t\}$. If there
is $\gt> t_0$ such that
$$f(\ors,p;\gt)> (|\ors|-m)(\gt-1) -m,$$
then $J(\ors|p)$ is finite.}

\medskip

We list only some results obtained by applying our Criterion
Theorem below. More examples including those when $l=1$ can be
found in Sec.~\ref{sec:gen}-\ref{sec:Jsd} and in the online
supplement \cite{onlineeg}.
\begin{thm} The $p$-divisible set $J(\ors|p)$ is finite if
\begin{quote}
\begin{description}
\item[$(1).$] $\ors=(1,1)$ and $p=3,7,13,31$, or $\ors=(1,1,1)$ and $p=3$,
or $\ors=(1,1,2)$ and $p=7$.

\item[$(2).$] $\ors=(4,3,5),\ (5,3,4)$  and $p=17$.

\item[$(3).$] $\ors=\{s\}^l$, $1\le l\le 20$, $s\ge 2$, and $p=2$.

\item[$(4).$] $\ors=(s,t)$, $s,t\le 20$, $t\ge 2$, and $p=2,3,5$.

\item[$(5).$] $\ors=(r,s,t)$, $r,s,t\le 10$, $t\ge 2$, and $p=2,3,5$.

\item[$(6).$]  $\ors=(q,r,s,t)$, $q,r,s,t\le 4$, $t\ge 2$, and $p=2,3,5$.
\end{description}
\end{quote}
Moreover, for $\ors$ in the last four cases we have
$J(\ors|2)=\{0\}$.
\end{thm}

\begin{conj}
For every $\ors$ the $2$-divisible set $J(\ors|2)=\{0\}$.
\end{conj}

In \cite[Conjecture B]{EL} Eswarathasan and Levine state that
there should be infinitely many primes $p$ (so called harmonic
primes) such that $J(1|p)=\{0,p-1,p^2-p,p^2-1\}$. Boyd \cite{Boyd}
further suggest $1/e$ as the expected density of such primes. For
any $\ors$ we extend this notion to define the {\em reserved
(divisibility) set} $RJ(\ors;x)$ of polynomials in $x$ with
rational coefficients. For any prime $p\ge |\ors|+3$  we have
$RJ(\ors;p)\subseteq J(\ors|p)$ and there are primes $p$ (called
{\em reserved primes} for $\ors$) such that equality holds. We
determine $RJ(\ors)$ for many types of $\ors$ in
Theorem~\ref{thm:RJ}. Furthermore, we argue by heuristics that the
following conjecture should be true.
\begin{conj}
If $l(\ors)=1$ and $\ors=s\ge 2$ then the proportion of primes $p$
with $J(s|p)=RJ(s)$ is $1/\sqrt{e}$. This proportion is equal to
$1/e$ for all other $\ors$.
\end{conj}
This conjecture is supported by very strong numerical and
theoretical evidence which we gather in \cite[Appendix II]{onlineeg} and in
Theorem~\ref{thm:RJ}.  It also generalizes Boyd's density conjecture of the
harmonic primes.

At the end of this paper we put forward some more
conjectures of $J(\ors|p)$ related to the distribution of
irregular primes.

This work was partially supported by NSF grant DMS0348258 and the Faculty
Development Fund from Eckerd College.

\section{A process to determine $J(\ors|p)$}\label{sec:gen}
For any positive integer $n$ let $n=p\tn+r$, where $\tn,r\in \N$
and $0\le r\le p-1$. For any $\ors\in \N^l$ define
\begin{equation*}
H^*(\ors;n)=\sum_{\substack {1\le k_1<\cdots<k_l\le n\\
p\nmid k_1k_2\cdots k_l}}  \frac{1}{k_1^{s_1}\cdots k_l^{s_l}} .
\end{equation*}
Then by a straightforward computation using the shuffle trick we
have: for any $s, l\in \N$
\begin{align}
H(s;n)=&H^*(s;n)+p^{-s}\cdot H(s;\tn),\label{*+tn} \\
H(\{s\}^l;n)=&\sum_{k=0}^l  p^{-ks}\cdot
    H(\{s\}^k;\tn)\cdot H^*(\{s\}^{l-k};n). \label{*+sN}
\end{align}
where $H(\{s\}^0;m)=H^*(\{s\}^0;m)=1$ for any integer $m$. We omit
the proofs of these formulas whose main ingredient is contained in
the proof of the main Criterion Theorem~\ref{thm:Jfinite} below.
Both of these formulas are generalizations of \cite[(2.2)]{EL} for
partial sums of harmonic series. They are the primary tools to
study Conjecture~\ref{conj:main} for homogeneous MHS.

For more general MHS we need a more complicated version of these
formulas. Fixing an arbitrary prime $p$ we define
$$G_0=\{0\}  \text{ and } G_t=\{n: p^{t-1}\le n< p^t\} \text{ for } t\in\N.$$
For any $\ors\in \N^l$ we set $J_t(\ors|p)=G_t\cap J(\ors|p).$

\begin{thm}\label{thm:Jfinite} {\em (\textit{Criterion Theorem})}
Let $l\ge 2$ be a positive integer and $p$ be a prime such that
$l\in G_{t_0}$. Let $\ors=(s_1,\cdots,s_l)\in \N^l$ and put
$m=\min\{s_i: 1\le i\le l\}$. For $t\in \N$ set
$f(\ors,p;t)=\min\{-v_p\big(H(\ors;n)\big): n\in G_t\}$. If there
is $\gt> t_0$ such that
$$f(\ors,p;\gt)> (|\ors|-m)(\gt-1) -m,$$
then $J(\ors|p)$ is finite.
\end{thm}
\begin{proof}
Let $n=p\tn+r \in G_{\gt+1}$. By definition we have
\begin{equation}\label{eqn:key}
H(\ors;n)=\sum_{\substack{\ga+\gb=l\\ \ga,\gb\ge 0} }
\sum_{\substack{1\le i_1<\dots<i_\ga \le n\\
1\le j_1<\dots<j_\gb \le n\\
\{i_1,\dots,i_\ga\}\cup\{j_1, \dots,j_\gb\}=\{1,\dots,n\} }}
\sum_{\substack
{1\le K_1<\dots<K_l \le n\\
p|K_{i_a}, 1\le a\le \ga \\
p\nmid K_{j_b}, 1\le b\le \gb }} \frac{1}{K_1^{s_1}\cdots
K_l^{s_l}}.
\end{equation}
Note that the terms corresponding to $\ga=l,\gb=0$ form the series
$A=H(\ors;\tn)/p^{|\ors|}\ne 0$ since $\gt>t_0$ so that $\tn\ge
p^{t_0}>l$. For all other terms with $\gb\ge 1$ we have the
natural bound
$$v_p( K_1^{s_1}\cdots  K_l^{s_l} )\le (|\ors|-m)\gt$$
since $K_j<p^{\gt+1}$ for all $j$ and one of them is prime to $p$.
Set $B=H(\ors;n)-A$ then obviously $v_p(B)\ge-(|\ors|-m)\gt.$
Since $\tn\in G_\gt$ by assumption of $\gt$ we know that
$$v_p(A)<-|\ors|-(|\ors|-m)(\gt-1)+m=-(|\ors|-m)\gt\le v_p(B).$$
By induction it's easy to see that for
all $t> \gt$ and $n\in G_t$ we have
$$v_p \big( H(\ors;n))<-|\ors|(t-\gt)
 -(|\ors|-m)(\gt-1)+m=-|\ors|(t-1)+m\gt<0.$$
This shows clearly that $J(\ors|p)$ is finite.
\end{proof}

\begin{cor} \label{cor:Jss}
Let $s$ and $l$ be two positive integers, $l\ge 2$, and $p$ be a
prime. Suppose $l\in G_{t_0}$. Then the $p$-divisible set
$J(\{s\}^l|p)$ is finite if there exists $\gt>t_0$ such that
$f_l(\gt):=f(\{s\}^l,p;\gt)>(l-1)s\gt-s$.
\end{cor}

\begin{prop} \label{prop:rs+rst}
Let $r,s,t$ be three positive integers. Let $\ors=(r,s)$ with
$1\le r\le 10$ and $2\le s\le 10$, or $\ors=(r,s,t)$ with $2\le
r,s,t\le 5$. Then there's always some prime $p\ge |\ors|+3$ such
that the $p$-divisible set $J(\ors|p)=RJ(\ors;p)$ is finite where
$RJ(\ors;p)$ is given in Theorem~\ref{thm:RJ}.
\end{prop}
\begin{proof}
The set $RJ(\ors;p)$ will be defined for general $\ors$ in
Definition~\ref{defn:RJ} and computed in Theorem~\ref{thm:RJ}. The
proof of the proposition follows from the Criterion
Theorem~\ref{thm:Jfinite} by computation. To save space we put the
details online \cite{onlineeg}.
\end{proof}

\section{Finiteness of $J(s|p)$}\label{sec:Js}
We now describe an approach to determine the $p$-divisible set $J(s|p)$
for any given
positive integer $s$ and odd prime $p$. This is essentially discovered
by Eswarathasan and Levine \cite{EL} and by Boyd \cite{Boyd},
independently. It follows quickly from Eq.~\eqref{*+tn} that
\begin{equation}\label{ItoJ}
n\in I(s|p) \mbox{ if and only if } \tn\in J(s|p^s).
\end{equation}
Therefore
\begin{equation}\label{oneway}
n\in J(s|p) \mbox{ implies }\tn\in J(s|p^s).
\end{equation}
It's also clear that
$$I(s|p)= p J(s|p^s) +R,\quad R=\{0, 1,\cdots, p-1\}.$$
\begin{rem}
The case when $s>1$ is very different from that of $s=1$
considered by previous authors in that the information of $I(s|p)$
is in general not enough to determine $J(s|p)$.
\end{rem}
To get an equivalent condition of Eq.~\eqref{oneway} we need
a partial generalization of \cite[Lemma~3.1]{EL}. Set the parity
function $\pari(m)=1$ if $m$ is odd and $\pari(m)=2$ if $m$ is even.
\begin{lem} \label{lem:pt}
Let $p$ be an odd prime and $s$ a positive integers. If $p-1\nmid
s,s+1$ then we have
\begin{equation}\label{Hpn}
H^*(s;pn)\equiv 0 \pmod{p^{\pari(s-1)}}.
\end{equation}
\end{lem}

\begin{proof}
By definition
\begin{equation}\label{H*}
H^*(s;pn)=\sum_{\substack {1\le k\le n\\ (p,k)=1}} \frac{1}{k^s}
=\sum_{m=0}^{n-1} \left(\sum_{mp< k< (m+1)p} \frac{1}{k^s}\right).
\end{equation}
The lemma follows from the fact that each inner sum in the parentheses
satisfies the congruence in Eq.~\eqref{Hpn} which can be proved by
the same argument as that in the proof of \cite[Lemma~2.2]{1stpart}
(when $n$ is odd the shorter proof suffices). It also follows from
\cite[Corollary~1]{sla}.
\end{proof}

\begin{prop}\label{prop:H*}
Let $s,l\in\N$ and  $p$ be an odd prime such that $p\ge l+2$ and
$p-1$ divides none of $ks$ and $ks+1$ for $k=1,\dots,l$. Then
$$H^*(\{s\}^l; pn)\equiv  0 \pmod{p^{\pari(ls-1)}}. $$
\end{prop}
The proof as well as the result itself is similar to that
of Theorem~\ref{th:p-1} so we leave the details to the interested reader.
The first step of induction is given as Lemma~\ref{lem:pt} above.
In fact, the proposition itself reduces to Theorem~\ref{th:p-1}
when $n=1$.
\begin{defn} For $n\in J(s|p^s)$ there is a unique
integer $\psi_s(s|p; n)\in [0,p-1]$ such that
\begin{equation}\label{psi}
\psi_s(s|p; n)\equiv \frac{1}{p^s} H(s;n) \pmod{p}.
\end{equation}
\end{defn}
\begin{lem}\label{lem:diff}
For $n=p\tn+r$, $0\le r<p$, we have
\begin{equation}\label{differen}
H^*(s;n)-H^*(s;p\tn)\equiv H(s;r) \pmod{p}.
\end{equation}
Furthermore, if $\tn\in J(s|p^s)$ then
\begin{equation}\label{difference}
H(s;n)\equiv H(s;r) + \psi_s(s|p; \tn) \pmod{p}
\end{equation}
\end{lem}
\begin{proof}
This follows from Eq.~\eqref{*+tn} and Lemma~\ref{lem:pt}.
Also see the proof of \cite[Lemma~3.2]{EL}.
\end{proof}
\begin{thm} \label{thm:equivalent}
Let $n=p\tn+r$, $0\le r<p$. Then $n\in J(s|p)$
if and only if
\begin{equation}\label{branching}
\tn\in J(s|p^s) \mbox{ and }  H(s;r)+\psi_s(\tn) \equiv 0 \pmod{p}.
\end{equation}
\end{thm}
\begin{proof} If $n\in J(s|p)$ then Eq.~\eqref{oneway}
implies that $\tn\in J(s|p^s)$. In addition, congruence~\eqref{branching}
follows immediately from Eq.~\eqref{difference}.
On the other hand, if Eq.~\eqref{branching} holds then
Eq.~\eqref{difference} implies that  $n\in J(s|p)$
and the proof is complete.
\end{proof}
We now use the above theorem to define a branching process
by using the sets $G_t$ which will compute $J(s|p)$ if it's finite.
\begin{prop} Let $s$ be a positive integer and $p$ an odd prime. Then
$J(s|p)=\cup_{t=0}^\infty J_t(s|p)$ where $J_t(s|p)$ can be
determined recursively by
$$J_{t+1}(s|p) =\{n=p\tn+r: \tn\in J_t(s|p^s), r\in R,
v_p\big(H(s;r)+\psi_s(\tn)\big)>0\}$$
for $t\in\N$. Here, as before, $R=\{0,1,\cdots, p-1\}.$
\end{prop}
The next corollary follows naturally.
\begin{cor} \label{cor:equivalent}
Let $s$ be a positive integer and $p$ an odd prime. Then $J(s|p)$ is
finite if and only if $J_t(s|p^s)=\emptyset$ for some $t\in\N$.
\end{cor}

An easy computation according to Corollary~\ref{cor:equivalent}
yields the following concrete result.
\begin{prop}\label{coneg}
 Let $s$ be a positive integer. Then $J(s|p)$ is
finite for primes $p=2,3,5,7$.
\end{prop}
\begin{proof}

(1) $p=2$. We claim that $J(s|2)=\{0\}$. We can prove that $2$ does not
divide $H(s;n)$ by induction on $n$. This is clear for $n=1$ and $n=2$
because $H(s;1)=1, H(s;2)=(1+2^s)/2^s$. Suppose $r\not\in J(s|2)$
for all $r\le n$ and $n\in J(s|2)$. If $n$ is odd then let
$H(s;n-1)=\frac{a}{2b}$ where $a$ is odd by inductive
assumption. Then
$$H(s;n)=\frac{a}{2b}+\frac{1}{n^s}= \frac{Na+2B}{\text{l.c.m}(2b,n^s)},$$
where $N=n^s/{\rm gcd}(n^s,b)$ and $B=b/{\rm gcd}(n^s,b)$.
Hence $Na+2B$ is odd because both $N$ and $a$ are odd, which is a
contradiction. If $n=2\tn$ then
$$H(s;2\tn)=\sum_{k=1}^\tn\left( \frac{1}{(2k-1)^s} + \frac{1}{(2k)^s}
\right) \equiv \tn+\frac{1}{2^s} H(s;\tn) \pmod{2}.$$ By inductive
assumption $2\nmid H(s;\tn)$ which implies that $2\nmid
H(s;2\tn).$ So $n$ can not belong to $J(s|2)$ either if $n$ is
even. This shows that $J(s|2)=\{0\}.$ In fact, it is not hard to
see that for $n\in G_t$, $t\ge 1$ we have
\begin{equation}\label{v2ofs1}
v_2\big(H(s;n)\big)= -(t-1)s.
\end{equation}

\medskip
For $3\le p\le 7$ Eswarathasan and Levine \cite{EL} have shown
that $J(1|p)$ are finite. We also know that when $s\le 4$ then
$J(s|p)$ are finite for these primes by explicit computation
\cite{onlineeg}. Assume $s\ge 4$. Then  by
Corollary~\ref{cor:equivalent} we only need to show that
$J_1(s|p^s)=\emptyset$. We need \cite[Corollary~2.7]{1stpart} which
implies that if $p\ge 3$ is a regular prime then
\begin{equation}\label{equ:sbig4}
H(s;p-1)\not \equiv 0 \pmod{p^s} \text{ for }s\ge 4.
\end{equation}

(2) $p=3$. Neither $H(s;1)=1$ nor $H(s;2)=1+1/2^s$ is divisible by $3^s$.
so $J_1(s|3^s)=\emptyset$.

(3) $p=5$. Neither $H(s;1)=1$ nor $H(s;2)=1+1/2^s$ is divisible by $5^s$.
Now
$$6^s H(s;3)=2^s+3^s+6^s \equiv 2^s+(-2)^s+1 \equiv
\begin{cases}
2\cdot 4^n+1  &\pmod{5}\quad\text{if $s=2n$},\\
1  &\pmod{5}\quad\text{if $s$ is odd}.
\end{cases}$$
So we always have $H(s;3)\not \equiv 0\pmod{5}$, i.e., $3 \not \in J_1(s|5^s)$.
Finally, Eq.~\eqref{equ:sbig4} implies that $4\not\in J_1(s|5^s)$ for $s\ge 4$
because 5 is a regular prime. Hence $J_1(s|5^s)=\emptyset$.

(4) $p=7$. Clearly $1,2\not \in J_1(s|7^s)$ and $6^s H(s;3)=2^s+3^s+6^s<7^s$
when $s\ge 4$. Now
$$H(s;4)=H(s;6)-\frac{1}{5^s}-\frac{1}{6^s}
\equiv (-1)^{s+1}\left(1+\frac{1}{2^s}\right) \pmod {7}.$$
Because $2^3\equiv 1 \pmod{7}$ we get
$$ 1+\frac{1}{2^s} \equiv \begin{cases}
2 &\pmod {7}\quad\text{ if } s\equiv 0\pmod{3},\\
3/2 &\pmod {7}\quad\text{ if } s\equiv 1\pmod{3},\\
3 &\pmod {7}\quad\text{ if } s\equiv -1\pmod{3}.
\end{cases}
$$
Therefore $4\not\in J_1(s|7^s)$. Similarly, $H(s;5)\equiv
(-1)^{s+1}\not\equiv 0 \pmod{7}.$ Finally, it follows
from Eq.~\eqref{equ:sbig4} that $6 \not \in J_1(s|7^s)$ for $s\ge 4$.
These show that $J_1(s|7^s)=\emptyset$ for all $s\ge 4$.
\end{proof}

\begin{rem}
The case $p=11$ is not so easy since $H(3;4)\equiv 0 \pmod{11}$ and
moreover, for any positive integer $e$ there is some $s<p^e(p-1)$ such
that $H(s;4)\equiv 0 \pmod{11^{e+1}}$.
\end{rem}

We also computed $J(s|p)$ for some other $s$ and $p$ (see
\cite{onlineeg}), which confirms the following
\begin{prop} \label{prop:eg}
Let $p$ be a prime such that $p\le 3001$. Then $J(s|p)$ is
finite for $2\le s\le 300$.
\end{prop}

\section{Finiteness of $J(\{s\}^l|p)$}\label{sec:Jsd}
In order to apply Criterion Theorem~\ref{thm:Jfinite} we set
$$f_l(t):=f(\{s\}^l,p;t)=\min\{-v_p(H(\{s\}^l ;n)): n\in G_t\}\quad
\forall t\ge 1.$$

We first look at the case $s\ge 2$.

\begin{lem}\label{lem:2s}
For all $s\ge 2$ we have $v_2(3^s+1)=\pari(s-1)$ which is 1 if $s$
is even and $2$ if $s$ is odd. In particular, we always have
$3^s+1\not\equiv 0\pmod{2^s}$.
\end{lem}
\begin{proof}
This is clear because
$$3^s+1 = \begin{cases}
9^n+1 \equiv 2        &\pmod{8} \quad\text{ if }s=2n,\\
3\cdot 9^n+1 \equiv 4  &\pmod{8}\quad\text{ if }s=2n+1.
\end{cases}$$
\end{proof}
\begin{prop} \label{prop:Jsd2}
Let $s\ge 2$ and $l\le 20$ be two positive integers. Then the
$2$-divisible set $J(\{s\}^l|2)$ is finite.
\end{prop}
\begin{proof}
When $l=1$ this is included in Proposition~\ref{coneg}. So we assume
$s,l\ge 2$. Then
$$H(s;2)= 1+ \frac{1}{2^s},\quad H(s,s;2)=\frac{1}{2^s},
\quad H(s;3)=  \frac{6^s+3^s+2^s}{6^s}.$$
Further, by Lemma~\ref{lem:2s} we know that
$$H(s,s;3)= \frac{1}{2^s} +\left(1+\frac{1}{2^s}\right)\frac{1}{3^s}
= \frac{3^s+1+2^s}{6^s}$$
has at least a factor 2 in the denominator. Therefore
we can take $\gt=2$ to get $f(\gt)\ge 1>(2s-s)(\gt-1)-s=0$.
So the condition in Corollary~\ref{cor:Jss}
is satisfied and consequently $J(s,s|2)=\{0\}$.

A detailed study of using Lemma~\ref{lem:2s} tells more. Let $t\ge
0$ and $n\in G_{t+2}$. Then by induction and Eq.~\eqref{*+sN}
we can easily show that
\begin{equation}\label{v2ofs2}
v_2\big(H(s,s;n)\big)=
\begin{cases}
-(2t+1)s &\quad\text{if } 2\cdot 2^t\le n<3\cdot 2^t, \\
\pari(s-1)-(2t+1)s &\quad\text{if } 3\cdot 2^t\le n<4\cdot 2^t.
\end{cases}
\end{equation}
Putting $l=3$ in the following equation
\begin{equation}\label{recc}
H(\{s\}^l;n)= H(\{s\}^l;n-1)+\frac{1}{n^s} H(\{s\}^{l-1};n-1),
\end{equation}
and applying induction on $t$ we can show that
\begin{equation}\label{v2ofs3}
v_2\big(H(s,s,s;n)\big)=
\begin{cases}
\pari(s-1)-3ts &\quad\text{if } 4\le 2n/2^t<5, \\
-3ts &\quad\text{if } 5\le 2n/2^t< 6 ,\\
-(3t+1)s &\quad\text{if } 6\le 2n/2^t<8.
\end{cases}
\end{equation}
So we get $J(s,s,s|2)=\{0\}$ when $s\ge 2$.

When $l\ge 4$ we can utilize Eq.~\eqref{recc} again. However, even in
the case $l=4$ it is very complicated already. Nevertheless the
idea is straightforward so we omit the details of the proof.
Suppose $s=2$ and $n\in G_{t+2}$ with $t\ge 1$ (note that
$H(\{s\}^4;n)=0$ for all $n\le 3$). Then we have
\begin{equation}
\label{v24}
v_2\big(H(\{2\}^4;n)\big)=
\begin{cases}
(1) \quad -2(4t-1)   &\quad\text{if } 32\le 16n/2^t<48, \\
(2) \quad -8t        &\quad\text{if } 48\le 16n/2^t<56, \\
(3) \quad -8t+\delta(t) &\quad\text{if } 56\le 16n/2^t <57 ,\\
(4) \quad -8t+7      &\quad\text{if } 57\le 16n/2^t<58,\\
(5) \quad -2(4t-2)   &\quad\text{if } 58\le 16n/2^t<60, \\
(6) \quad -2(4t-1)   &\quad\text{if } 60\le 16n/2^t<64.
\end{cases}
\end{equation}
Here if $t=1$ then (3)-(6) merge into (6); if $t=2$ then
$\delta(t)=5$ and (3)-(5) merge into (3); if $t=3$ then (3) and
(4) merge into (3); if $t\ge 3$ then $\delta(t)=6$. When $s=3$ and
$n\in G_{t+2}$ with $t\ge 1$ we have
\begin{equation}
\label{v34}
v_2\big(H(\{3\}^4;n)\big)=
\begin{cases}
(1) \quad -3(4t-1)   &\quad\text{if }
     8\le 4n/2^t<12, \\
(2) \quad -12t        &\quad\text{if }
    12\le 4n/2^t<14, \\
(3) \quad -3(4t-1)   &\quad\text{if }
    14\le 4n/2^t<15, \\
(4) \quad -3(4t-1)+1  &\quad\text{if }
    15\le 4n/2^t<16.
\end{cases}
\end{equation}
Here if $t=1$ then (3) and (4) merge into (4). When $s\ge 4$ we
have
\begin{equation}
\label{vs4}
v_2\big(H(\{s\}^4;n)\big)=
\begin{cases}
(1) \quad -s(4t-1)
    &\quad\text{if } 4\le 2n/2^t<6, \\
(2) \quad -4st
    &\quad\text{if } 6\le 2n/2^t<7, \\
(3) \quad -4st+2\pari(s-1)
    &\quad\text{if } 7\le 2n/2^t<8.
\end{cases}
\end{equation}
Equations \eqref{v24}-\eqref{vs4} imply that $J(\{s\}^4|2)=\{0\}$
for all $s\ge 2$.

Similar computation shows that when $l=5$ and $n\in G_{t+2}$ with
$t\ge 1$ we have
\begin{equation}
\label{v25}
v_2\big(H(\{2\}^5;n)\big)=
\begin{cases}
(1) \quad -s(5t-3)+2\pari(s-1)   &\quad\text{if } 16\le 8n/2^t<17, \\
(2) \quad -s(5t-3)+3   &\quad\text{if } 17\le 8n/2^t<18, \\
(3) \quad -s(5t-3)  &\quad\text{if }18\le 8n/2^t<20,\\
(4) \quad -s(5t-2)     &\quad\text{if } 20\le 8n/2^t<24,\\
(5) \quad -s(5t-1)   &\quad\text{if } 24\le 8n/2^t<28 , \\
(6) \quad -s(5t-1)+1
    &\quad\text{if } 28\le 8n/2^t<32.
\end{cases}
\end{equation}
Here if $t=1$ then (1)-(3) do not appear; if $s\ge 3$ then (1) and
(2)) merge into (1). This implies that $J(\{s\}^5|2)=\{0\}$ for
all $s\ge 2$.

As the length becomes longer (i.e. $l$ gets bigger) there are more
and more cases. The number of cases, denoted by $C(l)$, is
independent of $s$ when $s$ is large enough and tends to increase
with $l$ though not always. We compute the following
\begin{alignat*}{6}
\ &C(6)= 5,\quad  &\quad C(7)=7, \quad \ &\ C(8)=6,\quad  &\ C(9)= 8,\quad\ &\ C(10)= 8,\quad \\
\ &C(11)=11, &\quad C(12)=10,\ &\ C(13)=12,\ &\ C(15)=15,\ &\ C(16)=12,\\
\ &C(17)=15, &\quad  C(18)=14,\ &\ C(19)=18,\ &\ C(20)=15. \ & \
\end{alignat*}
After tedious verification we find $J(\{s\}^l|2)=\{0\}$ for all
$l\le 20$ and $s\ge 2$.
\end{proof}

For any given $l$ by similar method we should be able to determine
$J(\{s\}^l|2)$ for all $s\ge 2$. However, for odd primes $p$ we can
only extend this result to $J(\{s\}^l|p)$ for small $l$ and small $s$
with the aid of computers.

\begin{prop} \label{prop:Jsdbig}
Let $s$ and $l$ be two positive integers. Suppose $2\le s\le 10$
and $2\le l\le 10$. Then the $p$-divisible set $J(\{s\}^l|p)$ is
finite for the consecutive five primes immediately after $ls+2$.
Moreover there's always some prime $p$ such that
$J(\{s\}^l|p)=RJ(\{s\}^l;p)$ where
$$RJ(\{s\}^l;p)=\begin{cases}
\{0,p-1\} \quad &\text{if } 2\nmid s,\\
\{0,i+(p-1)/2,p-1: 0\le i\le l-1\} \quad &\text{if }2|s.
\end{cases}
$$
\end{prop}
\begin{proof}
The set $RJ(\ors;p)$ will be defined for general $\ors$ in
Definition~\ref{defn:RJ} and computed in Theorem~\ref{thm:RJ}. The
proof of the proposition follows from Corollary~\ref{cor:Jss} by
computer computation. To save space we put the details online
\cite{onlineeg}.
\end{proof}

In the rest of this section we turn to the case $s=1$. We may
assume $l\ge 2$  since the harmonic series has been handled by
\cite{EL} and \cite{Boyd}. According Corollary~\ref{cor:Jss} if we can
find $\gt$ large enough such that $f_l(\gt)\ge (l-1)(\gt-1)$ then
$J(1^l|p)$ is finite.

\begin{prop} \label{prop:J1l}
\begin{enumerate}
\item The $p$-divisible set $J(1^2|p)$ is finite if $p=3,7,13,31$.

\item  Let $s,t\le 20$ and $t\ge 2$. Then
the set $J(s,t|p)$ is finite for $p=2,3,5$.

\item  Let $r,s,t\le 10$ and $t\ge 2$. Then
the set $J(r,s,t|p)$ is finite for $p=2,3,5$.

\item  Let $q,r,s,t\le 4$ and $t\ge 2$. Then
the set $J(q,r,s,t|p)$ is finite for $p=2,3,5$.
\end{enumerate}
\end{prop}
\begin{proof}
We only need to find $\gt$ satisfying the condition of
Corollary~\ref{cor:Jss}.

(1) For each $\gt$ in the following we have
$f_2(\gt)=\gt-1$.

$p=3$. Take $\gt=6$. Then computation reveals that
$J(1,1|3)=\{0,5\}$. If $l=3$ then we take $\gt=10$. Then we have
$f_3(\gt)\ge 2(\gt-1)$. Note that in $G_9$ there is $n=17770$ such
that $v_3\big(H(1^3;n)\big)=-15$ so $f_3(9)=15$. By
Corollary~\ref{cor:Jss} and simple computation we see that
$J(1^3|3)=\{0,8\}$.

$p=7$. Take $\gt=4$. Then $J(1,1|7)=\{0,4,6,7,13\}$.

$p=13$. Take $\gt=4$. Then $J(1,1|13)=\{0,12,13,25\}$.

$p=31$. Take $\gt=4$. Then $J(1,1|31)=\{0,17,22,30,31,61\}$.

For the last three cases with $p=2,3,5$ we put the result of computation
online \cite{onlineeg}. For example, we can take
$\gt=10$ and show that $J(1,1,1|3)=\{0,8\}$.
\end{proof}
\begin{rem}
We could extend our results to larger $l$ and some other primes
$p$ but it would be very time consuming with our slow PCs. However,
even in the case $\ors=(1,1)$ similar process fails for $p=2$.
Computations suggest that $J(1,1|2)=\{0\}$,
$J(1,1|5)=\{0,4,5,9\}$, $J(1,1|11)=\{0,10,11,21\}$ and
$J(1,1|17)=\{0,11,13,16,17,33\}$. We will analyze the situation
for $p=2$ in detail in the next section.
\end{rem}

\section{Sequences related to $J(s,1|2)$}
One may wonder what goes wrong in Proposition~\ref{prop:J1l} if we
let $\ors=(1,1)$ and $p=2$. We will see that, amazingly, this
problem is closely related to some pseudo-random process.

Only in this section we adopt the shorthand $H_1(n):=H(1;n)$ and
$H_2(n):=H(1,1;n)$. Let's start with the first few partial sums of
$H_2(n)$ when $2\le n\le 14.$ Here $\sim$ means we only consider
the fractional part of the numbers.
\begin{align*}
\ &H_2(2)\sim\frac{1}{2}, \quad
H_2(3)\sim 1,\quad
H_2(4)\sim\frac{11}{24},\quad
H_2(5)\sim \frac{7}{8},\quad
H_2(6)\sim\frac{23}{90},\\
\ &H_2(7)\sim\frac{109}{180}, \quad
H_2(8)\sim\frac{9371}{10080},\quad
H_2(9)\sim\frac{467}{2016},\quad
H_2(10)\sim\frac{25933}{50400},\\
\ &H_2(11)\sim\frac{25933}{50400},\quad
H_2(12)\sim\frac{39353}{50400},\quad
H_2(13)\sim\frac{13501}{415800},\quad
H_2(14)\sim\frac{4027}{14850}.
\end{align*}
It looks like 2 never divides the numerator and moreover, the
$2$-powers in the denominators of $H_2(n)$ tend to increase with $n$,
though not always.
To proceed we need to know the $2$-divisibility of $H_1^*(n)$.
\begin{lem}\label{lem:H*}
Let $n$ be a positive integer. Then
\begin{equation*}
H_1^*(n)\equiv
\begin{cases}
0 \pmod{4} \quad&\text{ if } n\equiv 0,3 \pmod{4},\\
1 \pmod{4} \quad&\text{ if } n\equiv 1,2 \pmod{4}.
\end{cases}
\end{equation*} \end{lem}
\begin{proof}
If $n$ is even then obviously  $H_1^*(n)=H_1^*(n-1)$. So we
only need to consider $n\equiv 1, 3\pmod{4}$.

Set $\gd=1$ if $n=4l+1$ and $\gd=0$ if $n=4l-1$. Then
\begin{equation*}
H_1^*(n)=\frac{\gd}{4l+1}+\sum_{i=1}^{2l}\frac{1}{2i-1}
=\frac{\gd}{4l+1}+\sum_{i=1}^l\left(\frac{1}{2i-1}+\frac{1}{4l-2i+1}\right)
\equiv \gd \pmod{4}
\end{equation*}
as desired.
\end{proof}

\begin{rem} By working more carefully we can obtain the following
improvement of Lemma~\ref{lem:H*}:
if $n=2^dm$ or $n=2^dm-1$ where $m$ is odd and $d\ge 1$.
Then $v_2\big(H_1^*(n)\big)=2(d-1)$. However, the proof is complicated
and it is not needed in the rest of the paper so we leave the proof
of this general statement to the interested reader.
\end{rem}

The following result is exactly the reason why Corollary~\ref{cor:Jss}
cannot be applied to $J(1,1|2).$
\begin{prop}\label{prop:J11p2}
For any $t\ge 2$, there is a unique $n_t\in G_t$ such that
$v_2\big(H_2(n_t)\big)\ge 2-t$ whereas for all $n_t\ne n\in G_t$
we have $v_2\big(H_2(n)\big)\le 1-t$.  Therefore, for all positive
integers $n\not\in \{n_t\}_{t\ge 1}$ the numerator of $H_2(n)$ is
an odd integer.
\end{prop}
\begin{proof}
Note that $G_1=\{1\}$ and $G_2=\{2,3\}$. Thus $n_2=3$ because
$H_2(3)=1$. Assume that $t\ge 3$ and each $n_i$ has been found in
$G_i$ uniquely for $i\le t$. Let $n=2\tn+r\in G_{t+1}$ for $r=0$
or $1$. When $l=p=2$ and $s=1$ Eq.~\eqref{*+sN} becomes
\begin{equation}\label{equ:11p2}
H_2(n)= H_2^*(n)+\frac{1}{2} H_1(\tn)H_1^*(n)+\frac{1}{4}H_2(\tn).
\end{equation}
It's easy to show that $v_2\big(H_1(m)\big)=1-t$ for $m\in G_t$
by induction and the recursive relation $H_1(n)=H_1^*(n)+H_1(\tn)/2$.
If $\tn\ne n_t$ then we have $v_2\big(H_2(\tn)\big)\le 1-t$ and hence
$$v_2\big(H_2(n)\big)=\min\{v_2\big(H_1^*(n)\big)-t,
-1-t\}=-1-t<1-(t+1).$$
Suppose now $\tn=n_t$ and $n=2n_t+r_t$. We consider four possible
cases.

(i) If  $v_2\big(H_2(n_t)\big)=2-t$ and
$v_2\big(H_1^*(n)\big)\ge 1$ then $n\ne n_{t+1}$ because
$$v_2\big(H_2(n)\big)=\min\{v_2\big(H_1^*(n)\big)-t,
-t\}=-t=1-(t+1).$$

(ii) If  $v_2\big(H_2(n_t)\big)=2-t$ and
$v_2\big(H_1^*(n)\big)=0 $ then
$$v_2\left(\frac{1}{2} H_1(n_t)H_1^*(n)\right)=
v_2\left(\frac{1}{4}H_2(n_t)\right)=-t.$$
Hence $n=n_{t+1}$ because
$$v_2\big(H_2(n)\big)\ge 1-t=2-(t+1).$$

(iii) If $v_2\big(H_2(n_t)\big)\ge 3-t$ and $v_2\big(H_1^*(n)\big)=0 $
then $n\ne n_{t+1}$ because
$v_2\big(H_2(n)\big)=\min\{v_2\big(H_1^*(n)\big)-t,
v_2\big(H_2(n_t)\big)-2-t\}=-t=1-(t+1).$

(iv) If $v_2\big(H_2(n_t)\big)\ge 3-t$ and $v_2\big(H_1^*(n)\big)\ge 1 $
then $n=n_{t+1}$ because
$$v_2\big(H_2(n)\big)\ge \min\{v_2\big(H_1^*(n)\big)-t,
v_2\big(H_2(n_t)\big)-2-t\}\ge 1-t=2-(t+1).$$

\noindent
Now if $n_t=2l$ is even then by Lemma~\ref{lem:H*}

(1) $2n_t+1\equiv 1\pmod{4}$ and $v_2\big(H_1^*(2n_t+1)\big)=0$, and

(2) $2n_t\equiv 0\pmod{4}$ and $v_2\big(H_1^*(2n_t+1)\big)\ge 1$.

\noindent
If $n_t=2l+1$ is odd then by Lemma~\ref{lem:H*}

(3) $2n_t+1\equiv 3\pmod{4}$ and $v_2\big(H_1^*(2n_t+1)\big)\ge 1$ and

(4) $2n_t\equiv 2\pmod{4}$ and $v_2\big(H_1^*(2n_t+1)\big)=0$.

\noindent
Therefore, we have four situations to consider:

(a)  $n_t$ is even and  $v_2\big(H_2(n_t)\big)=2-t$. Then
$n_{t+1}=2n_t+1$ by (1) and (ii).

(b) $n_t$ is even and $v_2\big(H_2(n_t)\big)\ge 3-t$. Then
$n_{t+1}=2n_t$ by (2) and (iv).

(c) $n_t$ is odd  and  $v_2\big(H_2(n_t)\big)\ge 3-t$. Then
$n_{t+1}=2n_t+1$ by (3) and (iv).

(d) $n_t$ is odd  and  $v_2\big(H_2(n_t)\big)=2-t$. Then
$n_{t+1}=2n_t$ by (4) and (ii).

\noindent It follows that $n_{t+1}\in G_{t+1}$ is uniquely determined.
This finishes the proof of the proposition by induction.
\end{proof}

Denote the dyadic valuation $v_2\big(H_2(n_t)\big)$ by $-w_t$.
Then we have the following two interesting sequences:
\begin{align}
\ & \{n_t\}_{t\ge 2}=\left\{3,6,13,27,54,109,219,439,879,1759,3518,7037,14075,28151,
\right. \nonumber \\
\ & \qquad\qquad
\left. 56303,112606,225212,450424,900848,1801696,3603393,\dots \right\}\\
\ &\{w_t\}_{t\ge 2}=\{0,1,3,4,3,3,5,7,9,10,9,10,12,14,13,13,15,17,19,19,\dots\}
\end{align}
Set $r_1=r_2=1$ and
define $r_t=0$ or $1$ for $t\ge 3$ as determined in the proof of
Proposition~\ref{prop:J11p2} such that $n_{t+1}=2n_t+r_t$.
Then clearly $n_t$ can be written as
\begin{equation}\label{ntbin}
n_t=(r_1r_2\dots r_t)_2
\end{equation}
in binary system and apparently the sequence $\{n_t\}$ increases
very fast. Further, the occurrence of $r_t=0$ or $r_t=1$ does not
seem to have any predictable pattern and in fact we believe it is
related to some pseudo-random process. By this we mean the
following. First we of course conjecture that $\{w_t\}_{t\ge 3}$
is always bounded below by $1$ which is equivalent to say
$J(1,1|2)=\{0\}$. We also have proved that $w_t$ is bounded above
by $t-1$ and it is not hard to see that $w_t=t-1$ for infinitely
many $t$'s. It's also conceivable that $w_t$ are near $t-1$ most
of the time. However, we believe $w_t$ could move very far away
from $t-1$ for very large $t$. At the present stage, we could not
even determine whether the difference between $w_t$ and $t-1$ can
be arbitrarily large.

\begin{rem}
We put the two sequences $\{n_t\}_{t\ge 2}$ and $\{w_t\}_{t\ge 2}$
in Sloane's online database of integer sequences as $A079403(n)$
and $A079404(n)$, respectively. Shortly after Benoit Cloitre
emailed me a formula for the known terms of $\{n_t\}_{t\ge 2}$:
\begin{equation}\label{cconst}
n_t=\lfloor 2^{t-1}c \rfloor
\end{equation}
where $c=1.718232...$. Indeed, it's easy to see that
$$n_1=1,  n_t=2^{t-1} \prod_{k=1}^{t-1}\left(1+\frac{r_k}{2n_k}\right),
 \ \forall t\ge 2.$$
Further,
$$c=\lim_{t\to \infty} \prod_{k=1}^{t-1} \left(1+\frac{r_k}{2n_k}\right)$$
exists by comparison test. From Eq.~\eqref{ntbin} we get
$$c=(r_1.r_2r_3r_4\dots)_2=(1.101101111101111000001...)_2=1.718232... $$
Moreover, using binary system we see that the integral part of
$2^{t-1}c$ is exactly $n_t$, as desired.
\end{rem}

We can easily generalize Proposition~\ref{prop:J11p2} to the following.
\begin{prop}\label{prop:Js1p2}
For any $t\ge 2$, there is a unique $n_t\in G_t$ such that
$v_2\big(H(s,1;n_t)\big)\ge -s(t-1)+1$ whereas for all $n_t\ne
n\in G_t$ we have $v_2\big(H(s,1;n)\big)\le -s(t-1)$.  Therefore,
for every positive integer $n\not\in \{n_t\}_{t\ge 1}$ the
numerator of $H(s,1;n)$ is not divisible by $2$.
\end{prop}
\begin{proof} We can assume that $s\ge 2$ because of Proposition~\ref{prop:J11p2}.
The key to the proof is Eq.~\eqref{eqn:key} which yields that
$$H(s,1;n)=H^*(s,1;n)+U(s,1;n)+V(s,1;n)+\frac{1}{2^{s+1}}H(s,1;\tn)$$
where
$$U(s,1;n)=\frac{1}{2^{s}}\sum_{1\le 2k<l\le n, 2\nmid l}\frac{1}{k^sl},
 \quad V(s,1;n)=\frac{1}{2}\sum_{1\le k<pl\le n, 2\nmid k}\frac{1}{k^sl}.$$
Now it's easy to see that if $v_2\big(H(s,1;\tn)\big)\le  -s(t-1)$ then
$v_2\big(H(s,1;n)\big)\le -st$.

When $t=2$ we find $n_2=2$ always because
$$H(s,1;2)=\frac{1}{2}, \quad   H(s,1;3)=\frac{1}{2}+\frac{2^s+1}{3\cdot 2^s}.$$
Assume that $t\ge 3$ and $\tn$ is the unique $n_t\in G_t$
such that $v_2\big(H(s,1;\tn)\big)> -s(t-1)$. Then
$v_2\big(H(s;\tn)\big)=-s(t-1)$ by Eq.~\eqref{v2ofs1}. So we
can always uniquely choose $r_t$ so that for $n=2\tn+r_t\in G_{t+1}$
$$v_2\big(U(s,1;n)\big)=v_2\big(U(s,1;2\tn)+rH(s;\tn)/2^s\big)=-st$$
if $v_2\big(H(s,1;\tn)\big)=-s(t-1)+1.$   If
$v_2\big(H(s,1;\tn)\big)>-s(t-1)+1$ then we can uniquely choose $r_t$ so that
$$v_2\big(U(s,1;n)\big)\ge -st+1.$$
The upshot is for $n_t\in G_t$
there is a unique $n_{t+1}\in G_{t+1}$ satisfying the condition of
the proposition. This finishes the proof.
\end{proof}

In general, we cannot apply Criterion Theorem  to determine the
finiteness of $J((\ors,1)|2)$ for any $\ors\in \N^l$, because of
the existence of similar sequences. Moreover we believe $2$ never
divides the numerator of any multiple harmonic sum.
\begin{conj}
Let $l$ be an arbitrary positive integer and $\ors\in \N^l$. Then the $2$-divisible
set $J(\ors|2)=\{0\}$.
\end{conj}
We have verified this conjecture for all $\ors=(s,t)$ and
$\ors=(r,s,t)$ with $1\le r,s\le 10$ and $2\le t\le 10$, and for
all $\ors=(u,r,s,t)$ with $1\le u,r,s\le 4$ and $2\le t\le 4$. See
\cite{onlineeg}. The computation is very time-consuming, for
example when $\ors=(1,4,4,2)$ the Maple program runs more than 3.5
hours on my PC with Pentium 4 CPU 3.06GHZ and 512 MB RAM. The same
program in GP Pari runs a little faster. We put the program at the
end of our online supplement \cite{onlineeg}.

We believe that  among all possible $\ors$ and prime $p$ the cases
$(\ors,1)$ are the only ones that our Criterion Theorem fails (see
\cite{onlineeg}). Let me sketch a heuristic argument for this
belief for the case $\ors=1^l$ and $p\ge 3$.

By the recursive relation
$$H(1^l;n)=\sum_{j=0}^l  p^{-j}\cdot
    H(1^j;\tn)\cdot H^*(1^{l-j};n).$$
it is not hard to see that the size of $\gt$ we are looking
for in the Criterion Theorem depends on the length of the
sequences $\{n_t\}_{t>t_0}$ not satisfying the condition in the
theorem, where $n_t\in G_t$ and $n_{t+1}=pn_t+r_t$ for some $0\le
r_t<p$. If $n_t$ is already found then the existence of $n_{t+1}$
depends on $H^*(1;n)$ essentially, which we assume to distribute
among $(p+1)/2$ values modulo $p$ by the symmetric structure of
$J_1(1|p)$ (see Sec.~\ref{sec:J1}).
So $n_t$ produces two possible $n_{t+1}$ or no
$n_{t+1}$ with the same probability $q=(p-1)/(2p)$, and it produces
exactly one $n_{t+1}$ with probability $1/p$.

Let's assume that a certain cell reproduces itself according a
similar law as above, namely, it clones itself or dies in the next
generation with the same probability $q$, and with probability $1/p$ it stays alive
without reproduction. Let $p_k$ be the
probability that starting from $k$ cells in the beginning the
cells eventually all die out. We claim that $p_k=1$ for all $k$.
Indeed, it is not too hard to see that we only need to
show $p_1=1$. This follows from the Criticality Theorem for
Galton-Watson branching process (see \cite[Preface]{Gut} or
\cite[p.~7, Theorem~1]{AN})
because the average offspring is $2q+1/p=1$.

\section{The structure of $J_1(\{s\}^l |p)$}\label{sec:J1}
Set $J^0_t(\ors|m)=\{0\}\cup J_t(\ors|m)$ for any positive
integers $t$, $m$ and $\ors\in \N^l$. The next result is easy but
very useful in determining the structure of $J_1(s|p)$ since it
tells us essentially that $J^0_1(s|p)$ is symmetric about
$(p-1)/2$.
\begin{prop} \label{J1sym}
Let $p$ be an odd prime and $s\in \N$. Let $r\in \{1,\dots, p-2\}.$
Then $r\in J_1(s|p)$ if and only if $p-1-r\in J_1(s|p)$.
\end{prop}
\begin{proof}
If $p-1|s$ then $J_1(s|p)=\emptyset$ because
$H(s;r)\equiv r\pmod p$ for all $r\in \{1,\dots, p-2\}$.
If $p-1\nmid s$ then we have
$$H(s;r)=\sum_{k=1}^{r} \frac{1}{k^s}
=\sum_{k=p-r}^{p-1} \frac{1}{(p-k)^s}
\equiv (-1)^s \sum_{k=p-r}^{p-1} \frac{1}{k^s}\pmod{p}.$$
Subtracting $0\equiv (-1)^s \sum_{k=1}^{p-1} \frac{1}{k^s}\pmod{p}$
from the above we get the desired result.
\end{proof}

\begin{rem} We feel prompted to mention that the symmetry
of $ J^0_1(s|p)$ is not enjoyed by $J^0_1(s|p^2)$.
For instance, while $J^0_1(5|37)=\{0,6,9,12,18,24,27,30,36\}$
is symmetric about $18$ the set $J^0_1(5|37^2)=\{0,6,36\}$ is not.
\end{rem}

Now that we know $J^0_1(s|p)$ is symmetric we may wonder what
happens to the center $(p-1)/2$. When $s$ is odd the answer is
related to the irregularity of primes. 
\begin{prop} \label{halfway}
Let $p$ be an odd prime and $s$ be a positive integer such that
$p-1\nmid s$. Let $n$ be the unique positive
integer such that
$s\equiv n \pmod{p-1}$ and $2\le n\le p-2$. If
$n<p-3$ then we have
\begin{equation}\label{voron}
 H(s;(p-1)/2)\equiv
 \left\{
   \begin{array}{ll}
      \frac{2-2^{p}}{p}    \pmod{p}, & \hbox{if $n=1$;} \\
     \frac{2-2^{n}}{n}  B_{p-n} \pmod{p}, & \hbox{if $n>1$ is odd;} \\
      \frac{n(1-2^{n+1})}{2(n+1)}  pB_{p-n-1} \pmod{p^2}, & \hbox{if $n$ is even.}
   \end{array}
 \right. 
\end{equation}
Therefore,

(a) If $s$ is even then $(p-1)/2\in J_1(s|p)$. 

(b) If $s$ is odd and $(p,p-n)$ is an irregular pair then $(p-1)/2\in J_1(s|p)$.

(c) If $s$ is odd and $(p-1)/2\in J_1(s|p)$ and $2^s\not\equiv 2\pmod{p}$ then
$(p,p-n)$ is an irregular pair.

\noindent In particular, if $s\ge 3$ is odd and $p>2^s-2$ then
$(p-1)/2\in J_1(s|p)$ if and only if $(p,p-s)$ is an irregular
pair.
\end{prop}
\begin{proof} The congruence \eqref{voron} is essentially \cite[Cor.~5.2(b)]{Sun}.
The rest follows immediately.
\end{proof}

\begin{rem}\label{rem:halfway}
Proposition~\ref{halfway} tells us that when $s\equiv
1 \pmod{p}$ then $(p-1)/2\in J(s|p)$ if and only if $p$
is a Wieferich prime. There are only two known such primes:
$p=1093$ and $p=3511.$ If any others exist, they must be 
greater than $1.25×10^{15}$ according to \cite{K}.
\end{rem}

The above proposition says that if $s$ is even and $p-1\nmid s$ then
$(p-1)/2\in J_1(s|p)$. A natural question is when
$(p-1)/2\in J_1(s|p^2)$? The answer is given below.
\begin{cor} \label{cor:halfway}
Suppose $s$ is a positive integer such that $p-1\nmid s$.
Suppose $n$ is the integer between $2$ and $p-1$ such that
$n\equiv s \pmod{p-1}$ and $n<p-4$. For $n$ even 
$H(s;(p-1)/2)\equiv 0 \pmod {p^2}$
if and only if either $(p,p-n-1)$ is an irregular pair or
$p| 2^{n+1}-1$. If $n$ odd, $p|2^n-2$ and 
$(p,p-n)$ is an irregular pair then 
$H(s;(p-1)/2)\equiv 0 \pmod {p^2}$.  
\end{cor}
\begin{proof} This is clear for $n$ even. If $n$ is odd then
the corollary follows from \cite[Thm.~5.2(b)]{Sun}.
\end{proof}
\begin{rem}
For every positive even integer $s$ and every irregular prime $p> s+4$ up to
100,000, $p^2$ always exactly divides $H(s;(p-1)/2)$. Is this true in general?
The answer is no.  A calculation by Maple shows that for the 5952nd
irregular pair $(p,p-n-1)=(130811,52324)$ we have $n=78486$ and
$2^{n+1}\equiv 2 \pmod{p}$ and therefore $p^2|H(n+1;(p-1)/2)$ 
and $p^3|H(n;(p-1)/2)$ by \cite[Thm.~5.2(a),(b)]{Sun}.
The peculiarity of this pair was already noticed in \cite{EM}.
The next two such pairs are (599479,359568) (see \cite{EM2}), and
(2010401,1234960) (see \cite{BCEM}).
Note that apparently this problem is \emph{not} related to the problem of
$2^p\equiv 2\pmod {p^2}.$
\end{rem}

\begin{thm}\label{J12sd}
Let $s$ be a positive integer and $p>2ls+1$ be an odd prime. Then
$$\{p-1, j+(p-1)/2: j=0,1,\cdots,l-1\}\subset J_1\big(\{2s\}^l|p\big).$$
\end{thm}
\begin{proof}
Let $m=(p-1)/2$. By \cite[Lemma~2.12]{1stpart} there are integers
$c_\gl$ such that
\begin{equation}\label{clambda}
l! H\big(\{s\}^l;n\big)= \sum_{\gl\in P(l)} c_\gl H_\gl(s;n),
\end{equation}
where $P(l)$ is the set of partitions of $l$,
$H_{(\gl_1,\dots,\gl_r)}(s;n)=\prod_{j=1}^r H(\gl_j s;n)$ and
$c_{(l)}=(-1)^{l-1} (l-1)!.$ Plugging in $n=m$ we get $m\in
J_1\big(\{2s\}^l|p\big)$ by Proposition~\ref{halfway}. By definition
$H(\{s\}^l;q)=0$ for $q=1,\cdots,l-1$. When $q=1$ this implies
that $\sum_{\gl\in P(l)} c_\gl=0$ by Eq.~\eqref{clambda}. Hence
$m+1\in J_1\big(\{2s\}^l|p\big).$ Similarly, because $1$,
$1/2^{2s}$, $\cdots$, $1/(l-1)^{2s}$ are linearly independent when
regarded as a function of $s$, we see that for all independent
variables $x_1,\dots,x_j$, $j\le l-1$, we have
$$\sum_{\gl=(\gl_1,\dots, \gl_r)\in P(d)} c_\gl
 \prod_{j=1}^r (x_1^{\gl_1}+\cdots x_j^{\gl_j}) =0.$$
The theorem now follows from setting $x_j=1/(m+j)^{2s}$ for
$j=1,\cdots,l-1$.
\end{proof}

\begin{cor} \label{cor:2n2n}
Let $s$ be a positive integer and $p>4s+1$ be an odd prime. Then

(1) If $s$ is even then $(p-1)/2, (p+1)/2 \in J_1(s,s|p).$

(2) If $s$ is odd and $(p,p-s)$ is an irregular pair then $(p-1)/2
\in J_1(s,s|p). $ Further, if $s$ is odd, $2^s\not\equiv
2\pmod{p}$, and $(p-1)/2\in J_1(s,s|p)$, then $(p,p-s)$ is an
irregular pair. In particular,  if $s\ge 3$ is odd and $p>2^s-2$
then $(p-1)/2\in J_1(s,s|p)$ if and only if $(p,p-s)$ is an
irregular pair.
\end{cor}
\begin{proof}
Let $s$ and $p$ be the integers satisfying the conditions of the
corollary. When $s$ is even the corollary follows from
Theorem~\ref{J12sd}. If $s$ is odd then by Proposition~\ref{halfway} and
the shuffle relation we have
\begin{equation}\label{equ:2n2nm}
2H(s,s;(p-1)/2)\equiv H(s;(p-1)/2)^2-H(2s;(p-1)/2).
\end{equation}
The corollary follows immediately.
\end{proof}

\section{Reserved set of MHS}\label{sec:reserve}
In Conjecture B of \cite{EL} Eswarathasan and Levine state that
there should be infinitely many primes $p$ such that the divisible
set $J(1|p)=\{0,p-1,p^2-p,p^2-1\}$. Boyd \cite{Boyd} further
suggest $1/e$ as the expected density of such primes. The most
important steps are to elucidate the structure of $J_1(1|p)$ and
determine the relation between $J_t(1|p)$ and $J_{t+1}(1|p)$ for
$t>0$. We put forward some similar results and conjectures
concerning the divisible sets of general MHS in this last section.

\begin{defn}\label{defn:RJ}
For any $\ors\in \N^l$ there are finitely many function
$f_0(x)=0,f_1(x),\dots,f_r(x)\in \Q[x]$ such that for all primes
$p\ge wt(\ors)+3$
$$ f_0(p)<f_1(p)<\dots<f_r(p)  \text{ and }
f_0(p),f_1(p),\dots,f_r(p)  \in J(\ors|p).$$ We call the largest
$r$ the {\em reserved (divisibility) number} of MHS $H(\ors)$,
denoted by $\rho(\ors)$. We call the corresponding set
$\{f_0(x),\dots, f_{\rho(\ors)}(x)\}$ the {\em reserved
(divisibility) set} of $\zeta(\ors)$, denoted by
$RJ(\ors)=RJ(\ors;x)$. Its $t$-th segment is $RJ_t(\ors)=\{f(x)\in
RJ(\ors): f(p)\le p^t-1 \text{ for all prime } p\}$ for $t\ge 1$.
Note that $0\in RJ_t(\ors)$ for all $t\ge 0$. If
$J(\ors|p)=RJ(\ors;p)$ for some prime $p$ then is called a {\em
reserved prime} for MHS $H(\ors)$.
\end{defn}

For example, the reserved number of the harmonic series is 3,
the reserved set is $RJ(1)=\{0,x-1,x^2-x,x^2-1\}$, and 5 is
a reserved prime for the harmonic series because
$J(1|5)=\{0,4,20,24\}.$ We summarize
all known reserved sets in the following theorem.

\begin{thm}\label{thm:RJ}
Let $r,s,t,l\le 5$ be positive integers. Then
\begin{quote}
\begin{description}
\item[$(1).$] $RJ(1)=\{0,x-1,x^2-1,x^2-x\}.$

\item[$(2).$] If $1\ne l$ is odd  then $RJ_{10}(1^l)=\{0,x-1\}$.

\item[$(3).$] If $s\ge 3$ is odd then $RJ_{10}(\{s\}^l)=\{0,x-1\}$.

\item[$(4).$] If $s$ is even then $RJ_8(\{s\}^l) =\{0,i+(x-1)/2,x-1: 0\le
i\le l-1\}$.

\item[$(5).$] If $s\ge 3$ is odd then $RJ(s,1)=\{0,x-1,x\}$.

\item[$(6).$] If $s+t$ is odd then $RJ(s,t)=\{0\}$.

\item[$(7).$] If $s\ne t$, $s+t$ is even and $t\ne 1$ then
$RJ(s,t)=\{0,x-1\}$.

\item[$(8).$] \label{casestst} If $s\ne t$, $s+t$ is even and $t\ne 1$
then $RJ(s,t,s,t)=\{0,x-1\}$.

\item[$(9).$] If $r+s+t\ge 5$ is odd, $r,s,t\ge 2$, and $r\ne t$, then
$RJ(r,s,t)=\{0\}$.

\item[$(10).$] If $s$ is odd then $RJ(r,s,r)=\{0,x-1\}$.

\item[$(11).$] If $r+s+t\ge 6$ is even, $(r,s,t)\ne
(4,3,5),(5,3,4)$, and $r,s,t\ge 2$ are not all the same, then
$RJ(r,s,t)=\{0\}$. \label{H435conj}

\item[$(12).$] $RJ_{10}(2,1,1)=RJ(1,1,2) =RJ(4,3,5)=RJ(5,3,4)=\{0,x-1\}$.

\item[$(13).$] \label{m2n} If $s,t\ge 0$ and $s+t$ is even then
$RJ_{10}\big(1^s,2,1^t \big)=\{0,x-1\}$.

\item[$(14).$] \label{nnn} If $s\ge 1$ then
$RJ_{10}\big(1^s,2,1^s,2,1^s\big)=\{0,x-1\}$.

\item[$(15).$] \label{nn-1} If $s\ge 2$ is even then
$RJ_{10}\big(1^s,2,1^{s-1},2,1^{s+1}\big)=\{0,x-1\}$.

\item[$(16).$] \label{n-1n} If $s\ge 2$ is even then
$RJ_{10}\big(1^{s+1},2,1^{s-1},2,1^s\big)=\{0,x-1\}$.

\item[$(17).$] If $l$ is even then $RJ(1^l)=\{0,x-1,x,2x-1\}.$

\item[$(18).$] \label{case121} $RJ_{10}(1,2,1)=\{0,x-1,2x-1\}.$
\end{description}
\end{quote}
\end{thm}
\begin{rem}
We are sure that we can remove the subscript 10 in $RJ_{10}$
when more powerful computational tools are available to us.
\end{rem}
\begin{proof}
Even without the restriction of the bound 5, the inclusions
$RJ(\ors;p) \subseteq J(\ors|p)$ follow from Eq.~\eqref{clambda},
Theorem~\ref{th:p-1}, \cite[Theorems~2.14, 3.1, 3.5, 3.16, 3.18]{1stpart},
and Theorem~\ref{J12sd}, except in
case \eqref{casestst} and the last two cases.

For $\ors=1^l$, $l\ge 2$, Theorem~\ref{th:p-1} implies $p-1\in
J\big(1^l|p^{\pari(l-1)}\big)$ for all $l\ge 1$ and $p\ge l+3$. In
fact, by \cite[Theorems~2.8, 2.14]{1stpart} we have
\begin{equation}  \label{1dp-1}
H(1^l;p-1)\equiv \frac{H(l;p-1)}{(-1)^{l-1}l} \equiv
\begin{cases}
\displaystyle {-\frac{pB_{p-l-1}}{l+1} }  &\pmod{p^2} \quad\text{if } 2|l,\\
\displaystyle {-\frac{p^2(l+1)B_{p-l-2}}{2l+4} } &\pmod{p^3}
\quad\text{if } 2\nmid l.
 \end{cases}
\end{equation}
So if $p\ge l+3$  then
\begin{equation}  \label{1dp}
H(1^l;p)=H(1^l;p-1)+\frac{H(1^{l-1};p-1)}{p}\equiv
 \begin{cases}
\displaystyle {\frac{-p(l+2)B_{p-l-1}}{2(l+1)}} &\pmod{p^2} \quad\text{if } 2|l,\\
 \displaystyle {-\frac{B_{p-l}}{l} } &\pmod{p} \quad\text{if }2\nmid l.
 \end{cases}
\end{equation}
Further, setting $h_i=H(1^i;p-1)$, $h_0=1$ and $h_{-1}=0$ we have
\begin{multline}
 \label{onlyp2}
H(1^l;2p-1)=\sum_{i=0}l\ \sum_{1\le k_1<\dots <k_i\le
p <k_{i+1}<\dots<k_l<2p} \frac{1}{k_1k_2\dots k_l}  \\
\equiv \sum_{i=0}^l
 \Bigl(h_i+\frac{h_{i-1}}{p}\Bigr)
   \Bigl(h_{l-i}+(-1)^{l-i}pH(l-i+1; p-1) \Bigr)  \pmod{p^2}.
\end{multline}
Here we have used geometric series expansion inside $\Q_p$ such
that for any $m<p$
\begin{alignat*}{3}
\sum_{1\le i_1<\dots<i_m<p} \frac{1}{(p+i_1)\cdots (p+i_m)} \equiv
& h_m-p\sum_{i=1}^m H\big( (1^{i-1},2,1^{m-i});
p-1\big) &\pmod{p^3}\\
\equiv h_m+p((m+1)h_{m+1}-h_1h_m ) \equiv & h_m+(-1)^{m}pH(m+1;
p-1) &\pmod{p^3}
\end{alignat*}
by Eq.~\eqref{clambda}. It follows from equations
\eqref{clambda}, \eqref{onlyp2} and \cite[Theorem~3.1]{1stpart} that
\begin{equation*}
H(1^l;2p-1)\equiv 2h_l+\frac{1}{p}\sum_{i=0}^{l-1} h_i
h_{l-1-i}-(-1)^l H(l;p-1) \pmod{p^2}.
\end{equation*}
When $l$ is even we have
\begin{equation*}
H(1^l;2p-1)\equiv -\frac{(l+2)
H(l;p-1)}{l}+\frac{2H(l-1;p-1)}{p(l-1)} \equiv
 -2pB_{p-l-1} \pmod{p^2}.
\end{equation*}
So it's divisible by $p$. This shows that $2x-1$ belongs to the reserved
set $RJ(1^l)$ in the penultimate case
of the theorem. When $l=2n+1$ is odd $h_l\equiv 0
\pmod{p^2}$ and we get
\begin{equation*}
H(1^l;2p-1)\equiv \frac{1}{p}\sum_{i=0}^{n} h_{2i} h_{2n-2i}
\equiv  \frac{-H(2n;p-1)}{np} \equiv \frac{-2B_{p-l}}{l}\pmod{p}
\end{equation*}
which is rarely congruent to $0$. This explains why in case (2)
we can't have $2x-1$ in the reserved set $RJ(1^l)$ when $l$ is odd.

Finally let's turn to the last case of reserved set $RJ(1,2,1)$.
We have for any prime $p\ge 7$
\begin{equation*}
H(1,2,1;2p-1)= \sum_{1\le l<m<n<2p} \frac{1}{lm^2n}=A+B+C+D,
\end{equation*}
where
\begin{align*}
A=& \sum_{1\le l<m<n\le p} \frac{1}{lm^2n}= H(1,2,1;p-1)+\frac{1}{p}H(1,2;p-1), \\
B=& \sum_{1\le l<m\le p< n<2p} \frac{1}{lm^2n}=
 \left(H(1,2;p-1)+\frac{1}{p^2}H(1;p-1)\right)\sum_{1\le k<p}\frac{1}{p+k}, \\
C=& \sum_{1\le l\le p<m<n<2p} \frac{1}{lm^2n}=
\left(H(1;p-1)+\frac{1}{p}\right) \sum_{1\le m<n<p}\frac{1}{(p+m)^2(p+n)}, \\
D=& \sum_{p<l<m<n<2p} \frac{1}{lm^2n}=\sum_{1\le
l<m<n<p}\frac{1}{(p+l)(p+m)^2(p+n)}.
\end{align*}
Observe that $D\equiv H(1,2,1;p-1)\equiv 0 \pmod{p}$ by
\cite[Corollary~3.6]{1stpart} and $H(1;p-1)\equiv 0 \pmod{p^2}$ by
Wolstenholme's Theorem. By geometric series expansion we get
$$\sum_{1\le m<n<p}\frac{1}{(p+m)^2(p+n)}\equiv \sum_{1\le
m<n<p}\frac{(1-2p/m)(1-p/n)}{m^2n}  \pmod{p^2}.$$
Hence
\begin{multline*}
 H(1,2,1;2p-1)\equiv
\frac{1}{p}\big(H(1,2;p-1)+H(2,1;p-1)\big)\\
-2H(3,1;p-1)-H(2,2;p-1)\pmod{p}.
\end{multline*}
By Theorem~\ref{th:p-1} and \cite[Theorem~3.1]{1stpart} we have
$H(2,2;p-1)\equiv H(3,1;p-1)\equiv 0\pmod{p}.$ Further, from
shuffle relation we have
$$H(1,2;p-1)+H(2,1;p-1)=H(1;p-1)H(2;p-1)-H(3;p-1)\equiv 0\pmod{p^2}$$
by  Theorem~\ref{th:p-1}. This shows that $H(1,2,1;2p-1) \equiv
0\pmod{p}.$

To prove the theorem we now only need to demonstrate that
$RJ(\ors;p)=J(\ors|p)$ for some $p\ge |\ors|+3$ which can be done
through a case by case computation. We put this part of
verification online \cite{onlineeg}. In fact, much more data are
available in this supplement.
\end{proof}

Are there any other $\ors\in \N^l$ ($l\le 3$) besides those listed
in Theorem~\ref{thm:RJ} such that $RJ(\ors)\ne \{0\}$? In view of the
last conjecture of \cite{1stpart} we believe there are.

 From Theorem~\ref{thm:RJ} we see that to determine $RJ(\ors)$ we
often only need to study $RJ_1(\ors)$ because for all
non-homogeneous $\ors$ not equal to $(2r-1,1)$ or $(1,2,1)$, the
proportion of primes $p$ such that $J_t(\ors|p)=\emptyset$ for all
$t\ge 2$ is supposed to be positive. This implies that
$RJ(\ors)=RJ_1(\ors)$ for all such $\ors$. Precisely, we have
the following

\begin{conj}\label{conj:RJ1}
Suppose $\ors\in \N^l$ such that (i) $\ors=1$, or (ii)
$\ors=(1,2,1)$, or (iii) $\ors=(2r-1,1)$  for some $r\ge 1$, or
(iv) $\ors=1^{2l}$ for some $l\ge 1$. Then $RJ(\ors)=RJ_2(\ors).$
For all other $\ors$ we have $RJ(\ors)=RJ_1(\ors).$
\end{conj}

When $a\ge 2$ we assume that $H(\ors';ap-1)$ has random
distribution modulo $p^2$ (the case $a=2$ and $\ors'=(1,2,1)$ has
to be dealt with separately, but that's not hard). Thus the chance
that $p|H(\ors;ap)$ is less than $1/p^3$ for large $p$. This
implies that the probability of $J_2(\ors|p)=\emptyset$ is roughly
$(1-1/p^3)^{p^2}\to 1$ as $p\to \infty$.

\begin{defn} Let $\ors\in \N^l$. Define the {\em reserved density}
of the MHS $H(\ors)$ by
\begin{equation}
\text{density}(RJ(\ors); X)=\frac{\sharp\{\text{prime }p:
|\ors|+2< p<X,J(\ors|p)=RJ(\ors)\}} {\sharp\{\text{prime }p:
|\ors|+2< p<X\}}
\end{equation}
and the $m$th {\em reserved density}  by
\begin{equation}
\text{density}(RJ^m(\ors); X)=\frac{\sharp\{\text{prime }p:
|\ors|+1< p<X, \cup_{t=0}^m J_t(\ors|p)=RJ_m(\ors)\}}
{\sharp\{\text{prime }p: |\ors|+2< p<X\}}.
\end{equation}
\end{defn}

\begin{conj}\label{harmprimes}
Let $\ors\in \N^l$. Then
$$\operatorname{density}(RJ(\ors);\infty)=
\begin{cases}
1/\sqrt{e}, \quad &\text{ if $l=1,\ors\ge 2$},\\
1/e, &\text{ if $l=\ors=1$ or $l\ge 2$}.
\end{cases}
$$
\end{conj}

Note that we always have $RJ(\ors;p)\subseteq J(\ors|p).$ We have
put the data strongly supporting Conjecture~\ref{harmprimes} in
\cite[Appendix II]{onlineeg}. In fact, we have only computed the first or the
second reserved density because according to Conjecture~\ref{conj:RJ1}
this is enough to determine the reserved density in whole.

We now  provide a heuristic argument for Conjecture~\ref{harmprimes}.
Suppose $l=1$ and $\ors=s\ge 2$ first. Then by Proposition~\ref{J1sym}
we only consider $H(s;r)$ for $1\le r\le (p-5)/2+\pari(s-1)$
because for most $p$ the midpoint $(p-1)/2 \in J_1(s|p)$ if and
only if $s$ is even (see Proposition~\ref{halfway}). If we assume that
when $r$ varies the numbers $H(s;r)$ distribute randomly modulo
$p$ for any large fixed prime $p$ then the probability that
$J^0_1(s|p)=RJ_1(s;p)$ is $(1-1/p)^{(p-5)/2+\pari(s-1)}\to
1/\sqrt{e}$ as $p\to \infty$. By Conjecture~\ref{conj:RJ1} (we also
have a heuristic argument for it in this case) we see that the
probability that $J(s|p)=RJ(s;p)$ is $1/\sqrt{e}$.

\begin{rem} Although we cannot prove the random distribution of $H(s;r)$
for $1\le r\le (p-5)/2+\pari(s-1)$ modulo large prime $p$ ,
in a recent paper\cite{GFS},
Garaev et al. show that for any $\varepsilon > 0$,
the set $\{H(s;r)) : r< p^{1/2+\varepsilon}\}$ is
uniformly distributed modulo a sufficiently large $p$.
\end{rem}

Now we assume $l\ge 2$. In general $J^0_1(\ors|p)$ does not have
any symmetry so we see that the probability that
$J^0_1(\ors|p)=RJ_1(\ors;p)$ is $(1-1/p)^{p-\delta}\to 1/e$ as
$p\to \infty$, where $\delta=\sharp RJ_1(\ors).$ When $\ors$ does
not belong to the cases (i)-(iv) in Conjecture~\ref{conj:RJ1} we see
that the probability that $J(\ors|p)=RJ(\ors;p)$ is $1/e$ by
Conjecture~\ref{conj:RJ1}.

Finally let's deal with larger reserved sets when $l\ge 2$. By
Theorem~\ref{thm:RJ} we know that if $\ors=1^{2m}$ or $\ors=(1,2,1)$ or
$\ors=(2r-1,1)$ for some $r\ge 1$ then $RJ(\ors)=RJ_2(\ors)$. Let
$\ors=1^{2m}.$ Then for $p<n<2p-1$ by definition $H(1^{2m};n)$
is a sum of many product terms with $v_p$-value equal to either 0 or 1, i.e.,
the denominator has at most one $p$ factor. Assuming random distribution
of $pH(1^{2m};n)$ modulo $p^2$ we see that when
$p<n<2p-1$ the probability that $p$ divides $H(1^{2m};n)$ is
$1/p^2$. Thus the probability that none of $H(1^{2m};n)$ ($p<n<2p-1$)
is multiple of $p$ is equal to $(1-1/p^2)^{p-2}$.
Similar heuristic argument implies that
the probability that none of $H(1^{2m};n)$ ($2p\le n<3p$)
is multiple of $p$ is equal to $(1-1/p^3)^{p}$, and so on.
When  $l p\le n<(l+1)p$ and $2m\le l\le p$ the probability is
equal to $(1-1/p^{2m+1})^{p}$.
In general for $p^t\le n<p^{t+1}$ ($t\ge 1$) we can break it into
$p$ subintervals and estimate inside each of the subintervals.
It is easy to conclude that the probability that none of
$H(1^{2m};n)$ ($p^t\le n<p^{t+1}$)
is multiple of $p$ is at least $(1-1/p^{2m+2t-1})^{p^{t+1}}.$
Now
$$\sum_{t=1}^\infty \frac1{p^{2m+t-2}}= \frac1{p^{2m-2}(p-1)}\to 0\quad\text{as $p\to \infty$}. $$
We see that the probability that
$J(1^{2m}|p)=RJ(1^{2m};p)$ is
$$(1-1/p)^{p-2}\to \frac1e\quad\text{as $p\to \infty$}. $$
We omit the arguments for $\ors=(1,2,1)$ and $(2r-1,1)$
which are similar.

\medskip

We conclude our paper by some conjectures which concern
distributions of irregular primes in disguised forms.

\begin{conj}\label{conj:irrprimes} Let $r$, $s$ and $t$
be positive integers. Then
\begin{quote}
\begin{description}
\item[$(1).$]  If $s>1$ is odd then $J_1(s|p)=\{(p-1)/2,p-1\}$
for infinitely many primes $p$.

\item[$(2).$]  If $s>1$ is odd then
$J_1(s,s|p)=\{(p-1)/2,p-1\}$ for infinitely many primes $p$.

\item[$(3).$]  Let $s,t$ be positive integers. Suppose $s+t$ is odd.
Then $J_1(s,t|p)=\{p-1\}$ for infinitely many primes $p$.

\item[$(4).$] \label{conjd=3}
Let $r,s,t$ be positive integers such that $r+s+t$ is odd and $r\ne t$.
Then \makebox{$J_1(r,s,t|p)=\{p-1\}$} for infinitely many primes $p$.
\end{description}
\end{quote}
\end{conj}

Note that by various results of this paper and \cite{1stpart} an
affirmative answer to any part of our Conjecture~\ref{conj:irrprimes}
would imply that there are infinitely many irregular pairs
$(p,p-w)$ for any fixed odd number $w$ ($\ge 5$ in
case~\eqref{conjd=3}). Therefore, even if the sets of primes in
Conjecture~\ref{conj:irrprimes} are expected to be infinite they are
extremely sparse; very likely they have zero density.

\bigskip
\noindent
Department of Mathematics, Eckerd College, St. Petersburg,
Florida 33711, USA

\noindent
Max-Planck Institut f\"ur Mathematik, Vivatsgasse 7, 53111 Bonn, Germany

\noindent
zhaoj@eckerd.edu

\end{document}